\documentclass[12pt,a4paper]{amsart}
\usepackage{amssymb,amsmath,amsthm}

\usepackage{graphicx}

\usepackage{epsfig}

\long\def\onefigure#1#2{%  #1 picture,  #2  caption
\begin{figure*}[tbp]
\begin{center}
#1
\end{center}
\caption{#2}
\end{figure*}
} %end onefigure def

\newcommand{\lipefig}[2]  % labeled Ipe figure
{\onefigure{\mbox{\psfig{file=#1.eps}}}{\label{f:#1} #2} }

\begin{document}

\theoremstyle{plain}
\newtheorem{theorem}{Theorem}%[section]
\newtheorem{lemma}[theorem]{Lemma}
\newtheorem{prop}[theorem]{Proposition}
\newtheorem{corollary}[theorem]{Corollary}
\newtheorem{claim}[theorem]{Claim}

\newcommand{\sgn}{\textrm{sign}}
\newcommand{\tri}{\triangle}
\newcommand{\al}{\alpha}
\newcommand{\be}{\beta}
\newcommand{\ga}{\gamma}
\newcommand{\la}{\lambda}
\newcommand{\eps}{\varepsilon}
\newcommand{\N}{\mathbb{N}}
\newcommand{\z}{{\mathbb{Z}^2}}
\newcommand{\R}{\mathbb{R}}
\newcommand{\A}{\mathrm{A}}
\newcommand{\cc}{\mathcal{C}}
\newcommand{\xx}{\mathcal{A}}
\newcommand{\pp}{\mathcal{G}}
\newcommand{\dd}{\mathcal{D}}
\newcommand{\hh}{\mathcal{H}}
\newcommand{\dist}{\textrm{dist}}
\newcommand{\conv}{\textrm{conv}}
\newcommand{\intt}{\textrm{int}}
\newcommand{\sign}{\textrm{sign}}

\numberwithin{equation}{section}

\title{2013 unit vectors in the plane}

\author{Imre B\'ar\'any, Boris D. Ginzburg, and Victor S. Grinberg}

\keywords{unit vectors, norms, signed sum}

\subjclass[2010]{Primary 52A10, Secondary 15A39}

\begin{abstract} Given a norm in the plane and 2013 unit vectors in this norm, there is a signed sum of these vectors whose norm is at most one.
\end{abstract}

\maketitle
\bigskip
Let $B$ be the unit ball of a norm $\|.\|$ in $\R^d$, that is, $B$ is an 0-symmetric convex compact set with nonempty interior. Assume $V \subset B$ is a finite set. It is shown in \cite{BG} that, under these conditions, there are signs $\eps(v)\in \{-1,+1\}$ for every $v \in V$ such that $\sum _{v \in V}\eps(v)v \in dB$. That is, a suitable signed sum of the vectors in $V$ has norm at most $d$. This estimate is best possible: when $V=\{e_1,e_2,\dots,e_d\}$ and the norm is $\ell_1$, all signed sums have $\ell_1$ norm d.

In this short note we show that this result can be strengthened when $d=2$, $|V|=2013$ (or when $|V|$ is odd) and every $v \in V$ is a unit vector. So from now onwards we work in the plane $\R^2$.

\begin{theorem} Assume $V\subset \R^2$ consists of unit vectors in the norm $\|.\|$ and $|V|$ is odd. Then there are signs $\eps(v)\in \{-1,+1\}$ $(\forall v \in V)$ such that $\|\sum _{v \in V}\eps(v)v\| \le 1$.
\end{theorem}

This result is best possible (take the same unit vector $n$ times) and does not hold when $|V|$ is even.

\smallskip
Before the proof some remarks are in place here. Define the convex polygon $P=\conv\{\pm v:v\in V\}$. Then $P\subset B$, and $P$ is again the unit ball of a norm, $V$ is a set of unit vectors of this norm. Thus it suffices to prove the theorem only in this case.

\smallskip
A vector $v \in V$ can be replaced by $-v$ without changing the conditions and the statement. So we assume that $V=\{v_1,v_2,\dots,v_n\}$ and the vectors $v_1,v_2,\dots,v_n,-v_1,-v_2,\dots,-v_n$ come in this order on the boundary of $P$. Note that $n$ is odd. We prove the theorem in the following stronger form.

\begin{theorem} With this notation $\|v_1-v_2+v_3-\dots-v_{n-1}+v_n\| \le 1$.
\end{theorem}

{\bf Proof.} Note that this choice of signs is very symmetric as it corresponds to choosing every second vertex of $P$. So the vector $u=2(v_1-v_2+v_3-\dots-v_{n-1}+v_n)$ is the same (or its negative) when one starts with another vector instead of $v_1$. Define $a_i=v_{i+1}-v_i$ for $i=1,\dots,n-1$ and $a_n=-v_1-v_n$ and set $w=a_1-a_2+a_3-\dots+a_n$. It simply follows from the definition of $a_i$  that
\[
w=-2(v_1-v_2+v_3-\dots-v_{n-1}+v_n)=-u.
\]
Consequently $\|u\|=\|w\|$ and we have to show that $\|w\|\le 2$.

Consider the line $L$ in direction $w$ passing through the origin. It intersects the boundary of $P$ at points $b$ and $-b$. Because of symmetry we may assume, without loss of generality, that $b$ lies on the edge $[v_1,-v_n]$ of $P$. Then $w$ is just the sum of the projections onto $L$, in direction parallel with $[v_1,-v_n]$, of the edge vectors $a_1,-a_2,a_3,-a_4,\dots,a_n$. These projections do not overlap (apart from the endpoints), and cover exactly the segment $[-b,b]$ from $L$. Thus $\|w\|\le 2$, indeed. \hfill$\Box$

\smallskip
{\bf Remark.} There is another proof based on the following fact. $P$ is a zonotope defined by the vectors $a_1,\dots,a_n$, translated by the vector $v_1$. Here the zonotope defined by $a_1,\dots,a_n$ is simply
\[
Z=Z(a_1,\dots,a_n)=\left\{\sum_1^n\al_ia_i: 0\le \al_i\le 1 \;(\forall i)\right\}.
\]
The polygon $P=v_1+Z$ contains all sums of the form $v_1+a_{i_1}+\dots+a_{i_k}$ where $1\le i_1<i_2<\dots< i_k\le n$. In particular with $i_1=2,i_2=4,\dots,i_k=2k$
\[
v_1+a_2+a_4+\dots a_{2k}=v_1-v_2+v_3-\dots-v_{2k}+v_{2k+1} \in P.
\]
This immediately implies a strengthening of Theorem 1 (which also follows from Theorem 2).

\begin{theorem}\label{odd} Assume $V\subset \R^2$ consists of $n$ unit vectors in the norm $\|.\|$. Then there is an ordering $\{w_1,\dots,w_n\}$ of $V$, together with signs $\eps_i\in \{-1,+1\}$ $(\forall i)$ such that $\|\sum _1^k\eps_iw_i\| \le 1$ for every odd $k\in\{1,\dots,n\}$.
\end{theorem}

Of course, for the same ordering, $\|\sum _1^k\eps_iw_i\| \le 2$ for every $k\in\{1,\dots,n\}$. We mention that similar results are proved by Banaszczyk~\cite{Ba} in higher dimension for some particular norms.

\smallskip
In \cite{BG} the following theorem is proved. Given a norm $\|.\|$ with unit ball $B$ in $R^d$ and a sequence of vectors $v_1,\dots,v_n \in B$, there are signs $\eps_i\in \{-1,+1\}$ for all $i$ such that $\|\sum _1^k\eps_iw_i\| \le 2d-1$ for every $k\in\{1,\dots,n\}$. Theorem 1 implies that this result can be strengthened when the $v_i$s are unit vectors in $\R^2$ and $k$ is odd.

\begin{theorem} Assume $v_1,\dots,v_n\in \R^2$ is a sequence of unit vectors in the norm $\|.\|$. Then there are signs $\eps_i\in \{-1,+1\}$ for all $i$ such that $\|\sum _1^k\eps_iw_i\| \le 2$ for every odd $k\in\{1,\dots,n\}$.
\end{theorem}

The bound $2$ here is best possible as shown by the example of the max norm and the sequence $(-1,1/2),(1,1/2),(0,1),(-1,1),(1,1)$.

\smallskip
The {\bf proof} goes by induction on $k$. The case $k=1$ is trivial. For the induction step $k\to k+2$ let $s$ be the signed sum of the first $k$ vectors with $\|s\|\le 2$. There are vectors $u$ and $w$ (parallel with $s$) such that $s=u+w$, $\|u\|=1, \|w\|\le 1$. Applying Theorem 1 to $u,v_{k+1}$and $v_{k+2}$ we have signs $\eps(u),\eps_{k+1}$ and $\eps_{k+2}$ with $\|\eps(u)u+\eps_{k+1}v_{k+1}+\eps_{k+2}v_{k+2}\|\le 1$. Here we can clearly take $\eps(u)=+1$. Then
\[
\|s+\eps_{k+1}v_{k+1}+\eps_{k+2}v_{k+2}\|\le \|u+\eps_{k+1}v_{k+1}+\eps_{k+2}v_{k+2}\|+\|w\|\le 2
\]
finishing the proof. \hfill$\Box$

\bigskip
{\bf Acknowledgements.} Research of the first author was partially supported by ERC Advanced Research Grant no 267165 (DISCONV), and by Hungarian National Research Grant K 83767.

\vspace{1cm} {\sc Imre B\'ar\'any}
\\
  {\footnotesize R\'enyi Institute of Mathematics}\\[-1.5mm]
  {\footnotesize Hungarian Academy of Sciences}\\[-1.5mm]
  {\footnotesize PO Box 127, 1364 Budapest}\\[-1.5mm]
  {\footnotesize Hungary}\\[-1.5mm]
  {\footnotesize and}\\[-1.5mm]
  {\footnotesize Department of Mathematics}\\[-1.5mm]
  {\footnotesize University College London}\\[-1.5mm]
  {\footnotesize Gower Street, London WC1E 6BT}\\[-1.5mm]
  {\footnotesize England}\\[-1.5mm]
  {\footnotesize e-mail: {\tt barany@renyi.hu}}\\

\vspace{.4cm} {\sc Boris D. Ginzburg}
\\
{\footnotesize 1695 Betty Ct, Santa Clara, CA 95051}\\[-1.5mm]
  {\footnotesize USA}\\[-1.5mm]
  {\footnotesize e-mail: {\tt boris.d.ginzburg@gmail.com}}\\

\vspace{.4cm} {\sc Victor S. Grinberg}
\\
{\footnotesize 5628 Hempstead Rd, Apt 102, Pittsburgh, PA 15217}\\[-1.5mm]
  {\footnotesize USA}\\[-1.5mm]
  {\footnotesize e-mail: {\tt victor}\_{\;}grinberg@yahoo.com}\\

\end{document}